\documentclass{amsart}

\theoremstyle{plain}
\newtheorem{thm}{Theorem}
\newtheorem{prop}{Proposition}
\newtheorem{cor}{Corollary}
\newtheorem{lem}{Lemma}

\theoremstyle{definition}   
\newtheorem{de}{Definition}
\newtheorem{rk}{Remark}
\newtheorem{ex}{Example}
\newtheorem*{notation}{Notation}

\renewcommand{\iff}{if and only if}
\newcommand{\fs}{f-structure}
\newcommand{\ri}{Riemannian}
\newcommand{\ka}{K{\"a}hler}
\newcommand{\he}{Hermitian}
\newcommand{\mfd}{manifold}
\newcommand{\mo}{morphism}
\newcommand{\ha}{harmonic}
\newcommand{\ho}{holomorphic}
\newcommand{\hwc}{horizontally weakly conformal}
\newcommand{\ca}{Condition~(PHWC)}
\newcommand{\phm}{pseudo harmonic morphism}
\newcommand{\tm}{T^{\ast}M}
\newcommand{\tn}{T^{\ast}N}
\newcommand{\tc}{T^{\mathbb{C}}N}
\newcommand{\cn}{\mathbb{C}}
\newcommand{\rn}{\mathbb{R}}
\newcommand{\sn}{\mathbb{S}}
\newcommand{\co}{\circ}
\newcommand{\stm}{\scriptscriptstyle{\tm}}
\newcommand{\stn}{\scriptscriptstyle{\tn}}

\newcommand{\ssc}[1]{{\scriptscriptstyle{#1}}}

\DeclareMathOperator{\tr}{trace}

\def\operp{\bigcirc\kern-1.1em \bot }

\begin{document}

\begin{abstract}
We study a geometrical condition (PHWC) which is weaker than horizontal weak conformality. In particular, we show that {\ha} maps satisfying this condition, which will be called {\em {\phm}s}, include {\ha} {\mo}s and can be described as pulling back certain germs to certain other germs.
Finally, we construct a canonical {\fs} associated to every map satisfying (PHWC) and find conditions on this {\fs} to ensure the {\ha}ity of the map.
\end{abstract}

\title{Pseudo Harmonic Morphisms}
\author{E.~Loubeau}
\curraddr{Department of Pure Mathematics \\
            University of Leeds \\ 
           LS2 9JT Leeds \\
                 U.K.}
\address{(from 1st September 1996)\\ 
Universit{\'e} de Bretagne Occidentale \\
UFR Sciences et Techniques\\
Departement de Mathematiques\\
6, avenue Victor Le Gorgeu\\
BP 809\\
29285 BREST Cedex\\
France}
\email{pmt5el@amsta.leeds.ac.uk}
\maketitle

\section{Introduction}

Harmonic morphisms between {\ri} {\mfd}s are maps which pull back germs of {\ha} functions to germs of {\ha} functions. They can be characterised as {\hwc} {\ha} maps and the interplay between the analytical condition ({\ha}ity) and the geometrical one (horizontal weak conformality) is a rich source of properties.
\newline In this paper we shall be interested in generalising the geometrical condition of  horizontal weak conformality for maps into a {\he} {\mfd} to the condition of {\em Pseudo Horizontal Weak Conformality}
\begin{equation}
 [d\phi \circ (d\phi)^{\ast},J] = 0. \label{ast}
\end{equation}
Harmonic maps satisfying~\eqref{ast} shall be called {\em Pseudo Harmonic Morphisms}.
The origin of Condition~\eqref{ast} lies in the study of the stability of {\ha} maps. In particular, Burns, Burstall, de Bartolomeis and Rawnsley showed in~\cite{BurBurdeB89} that any stable {\ha} map into an irreducible {\he} symmetric space of compact type satisfies~\eqref{ast} and so, with our terminology, is a {\phm}. This result was then used by Chen~\cite{Che93A} to show that any stable {\ha} map from a compact {\ri} {\mfd} to the 2-sphere ${\sn}^{2}$ with its standard metric is a {\ha} {\mo}. 
\newline We first investigate the geometrical implications of~\eqref{ast} and explain how it generalises horizontal weak conformality. Then we show its independence of the complex structure on the target. 
\newline We then give a characterisation of {\em pseudo} {\ha} {\mo}s similar to that of {\ha} {\mo}s, as maps which pull back germs of pluri{\ha} functions to germs of {\ha} functions.
\newline Finally, to any map satisfying~\eqref{ast} we associate a generalised almost complex structure ({\fs}) and find conditions on this structure which force the map to be {\ha}.
\newline The author would like to wholeheartedly thank J.~C.~Wood for numerous discussions and his unfailing guidance.

\section[A Special Class...]{A Special Class of Harmonic Maps}

Let $(M,g)$ be a {\ri} {\mfd} of real dimension $m$, $(N,J,h)$ a {\he} {\mfd} of complex dimension $n$ and $\phi$ a smooth map between them.
\newline After complexifying the tangent bundle of $N$, we identify the tangent bundles of $M$ and $N$ with their duals, $\tm$ and $\tn$, in order to construct the map

$$ d\phi \circ (d\phi)^{\ast} : \tc \rightarrow \tc .$$

\begin{de}
We say that the map $\phi : (M,g) \rightarrow (N,J,h) $ is {\em Pseudo Horizontally Weakly Conformal\/} if
\begin{equation}
 [d\phi \circ (d\phi)^{\ast},J] = 0 . \tag{PHWC}
\end{equation}
\end{de}

\begin{lem}
{\ca} is equivalent to the condition that $d\phi \circ (d\phi)^{\ast}$ map the holomorphic tangent bundle $T^{(1,0)}N$ of $N$ onto itself.
\end{lem}

\begin{proof}
Let $v\in T^{(1,0)}N$, then, since $d\phi \circ (d\phi)^{\ast}$ is a complex linear map, {\ca} holds {\iff} 
\begin{equation*}
i(d\phi \circ (d\phi)^{\ast})(v) = J(d\phi \circ (d\phi)^{\ast})(v) \quad \forall v \in T^{(1,0)}N,
\end{equation*}
i.e. $(d\phi \circ (d\phi)^{\ast})(v)\in T^{(1,0)}N \quad \forall v \in T^{(1,0)}N$.
\end{proof}

\begin{lem} \label{Vlem2}
Set $V_{x} = d\phi_{x}^{\ast}({T^{(1,0)\ast}}_{\phi(x)}N)$ for $x\in M$. Then {\ca} is equivalent to $V_{x}$ being isotropic $\forall x\in M$.
\end{lem}

\begin{proof}
Let $v,w \in {T^{(1,0)\ast}}_{\phi(x)}N$, then
\begin{align*}
g_{x}(d\phi_{x}^{\ast}(v),d\phi_{x}^{\ast}(w) ) =& 
\, h_{\phi(x)}(v,d\phi_{x}\circ d\phi_{x}^{\ast}(w) ) ,
\end{align*}
so that
\begin{align*}
g_{x}(d\phi_{x}^{\ast}(v),d\phi_{x}^{\ast}(w) )=&\quad  0 \quad \text{for all $v,w \in {T^{(1,0)\ast}}_{\phi(x)}N$}
\end{align*}
{\iff} $d\phi\circ d\phi_{x}^{\ast}$ maps $T^{(1,0)\ast}N$ to itself i.e. {\iff} $\phi$ satisfies {\ca}.
\end{proof}

\begin{lem}
If we equip the manifold $(M,g)$ with the local real coordinates $(x^{i})_{i=1,\ldots,m}$ and the manifold $(N,J,h)$ with the local complex coordinates $(z^{\alpha})_{\alpha=1,\ldots,n}$ then {\ca} can be expressed as:
$$ \sum_{i,j=1}^{m}g^{ij}\frac{\partial\phi^{\alpha}}{\partial x^{i}}\frac{\partial\phi^{\beta}}{\partial x^{j}} = 0 \quad \forall\alpha,\beta =1,\ldots,n.$$
\end{lem}

\begin{proof}
We know that {\ca} is equivalent to $d\phi_{x}^{\ast}(T_{x}^{(1,0)\ast}N)$ being isotropic and this condition reads
\begin{align*}
g^{\ast}((d\phi)^{\ast}(dz^{\alpha}),(d\phi)^{\ast}(dz^{\beta})) = 0 \quad \forall \alpha,\beta = 1,\dots,n,
\end{align*}
where $g^{\ast}$ is the metric dual to $g$.
But the left-hand side of the last equation is equal to:
\begin{align*}
g^{\ast}\left(\frac{\partial\phi^{\alpha}}{\partial x^{i}}\, d x^{i},\frac{\partial\phi^{\beta}}{\partial x^{j}} \, d x^{j}\right) =& \\
 \sum_{i,j=1}^{m}\frac{\partial\phi^{\alpha}}{\partial x^{i}}\frac{\partial\phi^{\beta}}{\partial x^{j}}g^{\ast}(d
x^{i},d 
x^{j}) =& \sum_{i,j=1}^{m} g^{ij}\frac{\partial\phi^{\alpha}}{\partial x^{i}}\frac{\partial\phi^{\beta}}{\partial x^{j}}
\end{align*}
establishing the Lemma.
\end{proof}

\begin{rk} \label{rk3}
Let $(M,g)$ be a {\ri} {\mfd} and $(N,J,h)$ a {\he} {\mfd}. A map $\phi: (M,g)\rightarrow (N,J,h)$ is said to be {\em {\hwc}} if for each $p\in M$ either $d\phi_{p} = 0$ or $d\phi_{p}$ is conformal and surjective on the orthogonal complement of the kernel.. The conformal factor is then called the {\em dilation}. If $\phi$ is {\hwc} then, at all points, $d\phi \circ (d\phi)^{\ast} = \lambda^{2} Id$ (where $\lambda$ is the dilation of $\phi$ and $Id$ is the identity map on $\tc$). Therefore any {\hwc} map satisfies {\ca}, and, if $\dim_{\ssc\cn}N =1$, the two notions coincide. In particular any {\ha} morphism satisfies {\ca}.
\newline In local coordinates this phenomenon can be explained as follow. A map is {\hwc} {\iff} it satisfies
\begin{equation}
 g^{ij}\frac{\partial\phi^{A}}{\partial x^{i}}\frac{\partial\phi^{B}}{\partial x^{j}} = \lambda^{2} h^{AB} \quad\forall A,B=1,\overline{1},\ldots,n,\overline{n}.\label{hwc}
\end{equation}
If the metric $h$ is {\he} then $h^{AB} =0$ except for $A=\alpha$, $B=\overline{\beta}$ and $A=\overline{\alpha}$, $B=\beta$.
\newline If we decompose \eqref{hwc} into its $(1,1)$, $(2,0)$ and $(0,2)$ parts then 
{\ca} is clearly the $(2,0)$-part of \eqref{hwc}. Furthermore {\ca} is independent of the {\he} metric on $N$.
\newline
\end{rk}

\begin{de}
By a {\em local map} $\psi : N \to P$, we shall mean a map defined on an open subset of $N$.
\end{de}

\begin{prop}\label{lem6}
Let $(M,g)$ be a {\ri} {\mfd} (with local coordinates $(x^{i})$) and $N,P$ {\he} {\mfd}s (with respective systems of local coordinates $(w^{\alpha})$ and $(u^{a})$ ). 
Then the map $\phi: M\rightarrow N$ satisfies {\ca} {\iff} $\psi\circ\phi$ satisfies {\ca} for all local $\pm${\ho} maps $\psi : N\rightarrow P$.
\end{prop}

\begin{proof}
For any local {\ho} map $\psi$,
\begin{align}
g^{ij} \frac{\partial(\psi\circ\phi)^{a}}{\partial x^{i}} \frac{\partial(\psi\circ\phi)^{b}}{\partial x^{j}} = 
\frac{\partial\psi^{a}}{\partial w^{\alpha}}
\frac{\partial\psi^{b}}{\partial w^{\beta}}
g^{ij}
\frac{\partial\phi^{\alpha}}{\partial x^{i}}
\frac{\partial\phi^{\beta}}{\partial x^{j}}, & \label{c3}
\end{align}
with a similar formula for an anti{\ho} map. This clearly shows the ``if'' part, the converse being shown by choosing particular {\ho} maps.
\end{proof}

We can use {\ca} to define a class of maps which includes the notion of {\ha} morphism:

\begin{de}
A {\ha} map which is pseudo horizontally weakly conformal shall be called a {\em{Pseudo Harmonic Morphism}}.
\end{de}

\begin{prop} \label{lem4}
Let $\phi : (M,g) \rightarrow (N,J,h)$ be a smooth map from a {\ri} {\mfd} to a {\ka} {\mfd}.
\newline Then $\phi$ pulls back local complex-valued $\pm${\ho} functions on $N$ to local {\ha} functions on $M$ {\iff} $\phi$ is a {\phm}. In fact, in this case, $\phi$ pulls back local $\pm${\ho} functions to local {\ha} {\em morphisms}.
\end{prop}

\begin{proof}
All we need to do is to apply the chain rule for {\ha} maps and use {\ho} functions as test functions:
\newline Suppose that $\phi$ is a {\ha} map and let $f : N \to \cn$ be a local {\ho} function; then, with the usual notations:
\begin{align*}
 \tau(f\circ \phi) &= df(\tau(\phi)) + \tr\nabla df(d\phi,d\phi ) \\
&= \tr\nabla df(d\phi,d\phi ).
\end{align*}
As $(N,J,h)$ is {\ka}, in complex normal coordinates this can be written at a point as:
\begin{align*}
g^{ij}\frac{\partial^{2} f}{\partial z^{A}\partial z^{B}} \frac{\partial\phi^{A}}{\partial x^{i}}\frac{\partial\phi^{B}}{\partial x^{j}}=&\\
 g^{ij}\frac{\partial^{2} f}{\partial z^{\alpha}\partial z^{\beta}} \frac{\partial\phi^{\alpha}}{\partial x^{i}}\frac{\partial\phi^{\beta}}{\partial x^{j}} =& \quad\qquad\mbox{as $f$ is {\ho}}\\
=& \quad 0 \quad \mbox{if $\phi$ satisfies {\ca}.}
\end{align*}
Thus $f\circ \phi$ is {\ha}. Similarly for $f$ anti-holomorphic. 
\newline Conversely, if $\phi$ pulls back local $\pm${\ho} functions on $N$ to {\ha} functions on $M$ then 
$$ \frac{\partial f}{\partial z^{\gamma}} \tau^{\gamma}(\phi) + 
g^{ij}\frac{\partial^{2} f}{\partial z^{\alpha}\partial z^{\beta}} \frac{\partial\phi^{\alpha}}{\partial x^{i}}\frac{\partial\phi^{\beta}}{\partial x^{j}} = 0 $$
for any local {\ho} function $f$ on $N$. 
Choosing particular {\ho} functions shows that $\phi$ is {\ha} and satisfies {\ca}.

Let $f$ be {\ho}. Then by Proposition~\ref{lem6}, we know that $f\circ\phi$ satisfies {\ca} and because of Remark~\ref{rk3} this is equivalent to horizontal weak conformality so that $f\circ\phi$ is a {\ha} morphism.
\end{proof}

\begin{rk}
In the case that $(N,J,h)$ is only {\he} the proof fails as the expression for $\tr\nabla df(d\phi,d\phi )$ in local coordinates contains terms of the type
$$
g^{ij} \frac{\partial\phi^{\alpha}}{\partial x^{i}}\frac{\partial\phi^{\overline{\beta}}}{\partial x^{j}}
{}^{\ssc M}\Gamma^{\gamma}_{\alpha\overline{\beta}} \frac{\partial f}{\partial z^{\gamma}}
$$
which do not, in general, vanish. 
\end{rk}

After dealing with the case of {\ho} functions we turn our attention to the case of {\ho} maps.

\begin{prop} \label{lem5}
Let $(M,g)$ be a {\ri} {\mfd}. Let $(N,h)$ and $(P,k)$ be two {\ka} {\mfd}s. Then the map $\phi : (M,g)\rightarrow (N,h)$ is a {\phm} {\iff}  $\psi\circ\phi : (M,g)\rightarrow (P,k)$ is a local {\phm} for all local $\pm${\ho} maps $\psi : (N,h) \to (P,k)$..
\end{prop}

\begin{proof}
Given $p \in M$ and $\psi : N \to P$, equip the manifold $P$ with the complex normal local coordinates $(w^{m})_{m=1,\ldots,p}$ centred on $\psi(p)$ and $N$ with the complex normal local coordinates $(z^{\alpha})_{\alpha=1,\ldots,n}$ centred on $p$.
\newline Let $\phi : M \to N$ be a pseudo {\ha} {\mo} and $\psi : N \to P$ a local {\ho} map, then the map $\psi\circ\phi$ satisfies {\ca} because of Proposition~\ref{lem6}.
\newline We now show that $\psi\circ\phi$ is {\ha}:
\begin{align}
\tau^{m}(\psi\circ \phi) &= d\psi^{m}(\tau(\phi)) + \tr(\nabla d\psi)^{m}(d\phi,d\phi ) \notag\\
&= g^{ij}\frac{\partial^{2} \psi^{m}}{\partial z^{\alpha}\partial z^{\beta}} \frac{\partial\phi^{\alpha}}{\partial x^{i}}
\frac{\partial\phi^{\beta}}{\partial x^{j}} \qquad\text{since $\psi$ is {\ho}}\label{c2}\\
& =  \quad 0 \notag
\end{align}
so $\psi\circ\phi$ is {\ha}.
\newline Conversely, if $\psi\circ\phi$ is a {\phm} for any local $\pm${\ho} map $\psi$ then, again because of Proposition~\ref{lem6}, $\phi$ satisfies {\ca}.
\newline Further, as $\psi\circ \phi$ is {\ha}, we have
\begin{equation*}
\tau^{m}(\psi\circ \phi) = d\psi^{m} (\tau(\phi)) + \tr{ \nabla d\psi \left(d\phi, d\phi \right) }
=\, 0 .
\end{equation*}
As $\psi$ is $\pm${\ho} and $\phi$ satisfies {\ca} the second term is zero (cf.~\eqref{c2}) and, by taking particular {\ho} maps, we see that $\phi$ is {\ha}.
\end{proof}

\begin{rk}
The ``if'' part of Proposition~\ref{lem4} is a special case of the ``only if'' part of Proposition~\ref{lem5}.
\end{rk}

\begin{de} \label{defi4}
A smooth map from a complex {\mfd} $N^{n}$ to a {\ri} {\mfd} $P^{r}$ satisfying:
\begin{equation}
\nabla^{(0,1)}d'\phi = 0 \label{IVequaP}
\end{equation}       
is called a {\em{pluri{\ha}  map}} (\cite{Uda88A}). Here
$$\left( \nabla^{(0,1)}d'\phi \right) (\overline{Z},W) = 
{}^{\phi}\nabla_{\overline{Z}} \left( d'\phi(W) \right) - 
d'\phi \left( \overline{\partial}_{\overline{Z}} W \right) \qquad \forall Z,W\in T^{(1,0)}N^{n},$$
where $\overline{\partial}$ is the $\overline{\partial}$-operator of the holomorphic bundle $T^{(1,0)}N^{n}$.
A map is a  pluri{\ha} map if and only if its restriction to any holomorphic curve is harmonic (see~\cite[Prop. 1.1]{OhnVal90}).
\newline We remark that:
$$\nabla^{(0,1)}d'\phi = 0 \Leftrightarrow \nabla^{(1,0)}d''\phi =0 . $$
The equation for pluriharmonicity of a map is
\begin{align} 
&\frac{\partial^{2} \phi^{A}}{\partial z^{\alpha} \partial z^{\overline{\beta}}}
+ {}^{\ssc X}\Gamma^{\ssc{A}}_{\ssc{JL}} 
\frac{\partial\phi^{\ssc{J}}}{\partial z^{\alpha}} \frac{\partial\phi^{\ssc{L}}}{\partial z^{\overline\beta}} = \, 0 \tag{PH}\label{equaph}\\
&\forall A\in \{1,\dots ,r \} \mbox{ and } \forall \alpha, \beta\in \{1,\dots ,n\}.\notag
\end{align}

\end{de}

\begin{rk}
Just like ``being {\ha} and satisfying {\ca}'' is independent of any choice of metric on $N$ (see Proposition~\ref{lem5} with $(P,k) = (N,h)$), pluri{\ha} maps are
harmonic with respect to any {\ka} structure on the manifold $N$.
\end{rk}

\begin{prop} \label{V6}
Let $\phi$ be a smooth map from a {\ri} {\mfd} $(M,g)$ to a {\ka} {\mfd} $(N,h)$.
Then $\phi$ pulls back local pluri{\ha} functions on $N$ to local {\ha} functions on $M$ {\iff} $\phi$ is a {\phm}.
\end{prop}

\begin{proof}
As a special case of Definition~\ref{defi4}, a pluri{\ha} function on a complex {\mfd} $N$ is a map $f:N\rightarrow \rn$ which satisfies
$$\nabla^{(0,1)}d' f = 0.$$
In complex normal local coordinates this condition is (cf. \cite[p. 631]{OhnUda90})
$$ \frac{\partial^{2} f}{\partial z^{\alpha}\partial z^{\overline\beta}} =0 \qquad\forall \alpha,\beta =1,\ldots,n.$$
All we need to do is to apply the chain rule for {\ha} maps:
\newline Suppose that $\phi$ is a {\phm} and let $f : N \rightarrow \rn$ be a local pluri{\ha} function. Then, in normal local coordinates:
\begin{align*}
 \tau(f\circ \phi) &= g^{ij}\frac{\partial^{2} f}{\partial z^{A}\partial z^{B}} \frac{\partial\phi^{A}}{\partial x^{i}}\frac{\partial\phi^{B}}{\partial 
x^{j}}\\
&= g^{ij}\frac{\partial^{2} f}{\partial z^{\alpha}\partial z^{\beta}} \frac{\partial\phi^{\alpha}}{\partial x^{i}}\frac{\partial\phi^{\beta}}{\partial x^{j}} + 
g^{ij}\frac{\partial^{2} f}{\partial z^{\overline\alpha}\partial z^{\overline\beta}} \frac{\partial\phi^{\overline\alpha}}{\partial x^{i}}\frac{\partial\phi^{\overline\beta}}{\partial x^{j}}
 \quad\qquad\mbox{as $f$ is pluri{\ha}}\\
=& \quad 0 \quad \mbox{as $\phi$ satisfies {\ca}.}
\end{align*}

To prove the converse, first observe that a complex-valued function is pluri{\ha} {\iff} its real and imaginary parts are real-valued pluri{\ha} functions, and that $\pm${\ho} functions on $N$ are pluri{\ha}. Therefore by taking the real (or imaginary) part of a local $\pm${\ho} function on $N$ with prescribed first derivatives at a point $p\in N$,  we can produce local real-valued pluri{\ha} functions on $N$ with prescribed derivatives at the point $p$.
\newline If $\phi$ pulls back local pluri{\ha} functions on $N$ to {\ha} functions on $M$ then 
$$ \frac{\partial f}{\partial z^{\gamma}} \tau^{\gamma}(\phi) + 
g^{ij}\frac{\partial^{2} f}{\partial z^{\alpha}\partial z^{\beta}} \frac{\partial\phi^{\alpha}}{\partial x^{i}}\frac{\partial\phi^{\beta}}{\partial x^{j}} = 0 $$
for any local pluri{\ha} function $f$ on $N$.
\newline We then choose pluri{\ha} functions with prescribed first derivatives at a point, by the process described above, and proceed as in Proposition~\ref{lem4}.
\end{proof}

Similarly,
\begin{cor}
Let $N$ and $P$ be {\ka} {\mfd}s.
Then a smooth map $\phi : M \to N$ pulls back local pluri{\ha} maps from $N$ to $P$ to {\ha} maps on $M$ {\iff} it is a {\phm}.
\end{cor}

We close this section with examples which show that {\ca} and horizontal weak conformality  are not equivalent.

\begin{ex}
Let $\phi : \rn^{2} \rightarrow {\cn}^{\,3}$ be defined by:
$$\phi(x_{1},x_{2}) = 
\left( \phi^{1}(x_{1},x_{2}),\phi^{2}(x_{1},x_{2}),\phi^{3}(x_{1},x_{2}) \right)
$$
with
$$\phi^{1}(x_{1},x_{2})=\phi^{2}(x_{1},x_{2})=\phi^{3}(x_{1},x_{2})= x_{1} + ix_{2}.$$
It is easy to see that $\phi$ is {\ha} and that it cannot be {\hwc} because it is not constant and the dimension of the target is greater than the dimension of the domain.
\newline However
\begin{align*}
\sum_{i=1,2} \frac{\partial\phi^{\alpha}}{\partial x_{i}}
\frac{\partial\phi^{\beta}}{\partial x_{i}} = 1^{2} + i^{2} = 0 , \quad \alpha,\beta = 1,2,3.
\end{align*}
So $\phi$ is {\ha} and satisfies {\ca} but is not {\hwc} (it is, in fact, an immersion).
\end{ex}
\begin{ex}
Let $\phi : \rn^{4} \rightarrow {\cn}^{\,2}$ be defined by:
$$\phi(x_{1},x_{2},x_{3},x_{4}) = 
\left( i(x_{1} + x_{2}) + x_{3} + x_{4},  i(x_{1} + x_{2}) + x_{3} + x_{4} \right).$$
It is easy to see that 
$$\left( \ker{d\phi} \right)^{\bot} = \left\{ \left( X_{1}, X_{2}, X_{3}, X_{4} \right) \vert \, X_{1}=X_{2} , X_{3} = X_{4} \right\}.$$
\newline If $X,Y \in \left( \ker{d\phi} \right)^{\bot}$ then
$$ \left< X , Y \right> = 2\left( X_{1}Y_{1} + X_{3} Y_{3} \right)$$
while 
$$\left< d\phi(X), d\phi(Y) \right> = 8 \left( X_{1}Y_{1} + X_{3} Y_{3} + i \left( X_{1}Y_{3} - X_{3}Y_{1} \right) \right), $$
so that $\phi$ cannot be {\hwc}.
\newline However $\phi$ is {\ha} since linear and it is easy to verify that $\sum_{i=1}^{4} \left( \frac{\partial \phi^{\alpha}}{\partial x_{i}} \right)^{2} = 0$, i.e. $\phi$ is a {\phm}. 
\end{ex}

\section{{f}-Structures and (PHWC) Maps}

\begin{de}
An {\fs} on a {\ri} {\mfd} $(M,g)$ is a (smooth) skew-symmetric section $F$ of End$(TM)$ such that:
\begin{equation}
 F^{3} + F = 0 . \label{F}\tag{F}
\end{equation}
\end{de}

The definition implies that $F$ has three possible eigenvalues on the complexification of $TM$: $+i$, $-i$ and $0$. The corresponding eigenspaces $T^{+}M$, $T^{-}M$ and $T^{0}M$ are orthogonal with respect to the {\he} metric $h(X,Y)=g(X,\overline{Y})$ on $T^{\ssc\cn}M$. An almost complex structure on a {\mfd} is an {\fs} with trivial $0$-eigenspace.

f-Structures were first considered by Yano~\cite{Yano63} (cf.~\cite{IshiharaYano64} as well as \cite[Chap.~VII]{KonYano84}). They include almost complex structures and almost contact structures.

\begin{prop} \cite[Prop.~2.2]{Raw85} \label{Vbijec}
Let $(M,g)$ be a {\ri} {\mfd}.
There is a bijection between the {\fs}s $F$ of rank $2k$ on $TM$ and the $g$-isotropic subbundles $T^{+}M$ of $T^{\ssc\cn}M$ of rank $k$, given by $T^{+}M = \ker(F - i)$.
\end{prop}

\begin{de}\cite{IshiharaYano64}
An {\fs} $F$ of rank $2m$ is said to be {\em integrable} if there exists a system of local coordinates {\em adapted} to $F$, i.e. in which $F$ has the constant components

$$ F = 
\begin{pmatrix}
0 & -\mathrm{Id}_{m} & 0 \\
\mathrm{Id}_{m} & 0 & 0 \\
0 & 0 & 0
\end{pmatrix},
$$
where $\mathrm{Id}_{m}$ is the $m\times m$ identity matrix.
\newline An {\fs} $F$ is said to be {\em parallel} if
$$ (\nabla F) (X,Y) =0 \quad \forall X,Y \in T^{\ssc\cn}M.$$
As in the case of complex structures, parallelism implies integrability.
\end{de}

\begin{thm} \cite{IshiharaYano64} \label{Vintegrability}
Let $F$ be an {\fs}. 
Define the {\em Nijenhuis tensor} $N(X,Y)$ of $F$ by
$$N(X,Y) = [FX,FY] - F[FX,Y] - F[X,FY] + F^{2}[X,Y] .$$
Then the {\fs} $F$ is integrable {\iff} 
$$ N(X,Y) = 0 $$
for any two vector fields $X$ and $Y$.
\end{thm}

\begin{de}
A map $\phi$ between two {\ri} {\mfd}s $(M,g)$ and $(N,h)$, each carrying an {\fs} $F^{M}$ and $F^{N}$, shall be called {\em f-{\ho}} if $d\phi$ intertwines the {\fs}s, i.e.
$$d\phi \co F^{M} = F^{N} \co d\phi.$$
\end{de}

\begin{rk}
\begin{enumerate}
\item This definition implies that
\begin{align*}
(d\phi)(T^{+}M) \subseteq T^{+}N ,\,
(d\phi)(T^{-}M) \subseteq T^{-}N ,\,
(d\phi)(T^{0}M) \subseteq T^{0}N.
\end{align*}
\item This choice of terminology differs slightly from Rawnsley's in \cite{Raw85} where he prefers to call such maps ``{\mo}s of f-{\mfd}s'', reserving the term ``f-{\ho}'' for the cases where the domain carries an almost complex structure.
\end{enumerate}
\end{rk}

An application of {\fs}s to the theory of {\ha} maps was given by Rawnsley in the following  generalisation of a result of Lichnerowicz:

\begin{prop} \cite{Raw85}
Let $(M,g,J)$ be a cosymplectic {\mfd} and let $(N,h,F)$ be a {\ri} {\mfd} with an {\fs} which satisfies:
$$ \nabla_{X} C^{\infty}(T^{+}N)  \subseteq C^{\infty}(T^{+}N) \qquad \forall X\in C^{\infty}(T^{-}N).$$
Then every f-{\ho} map $\phi :M\to N$ is {\ha}.
\end{prop}

\begin{de}
We shall say that a complex vector field $v \in C(T^{\ssc\cn}M)$ is of type $+$ (respectively $-$, $0$) if $v\in C(T^{+}M)$ (respectively $v\in C(T^{-}M)$, $v\in C(T^{0}M)$). Similarly for complex 1-forms $\theta \in C(T^{{\ssc\cn}\ast}M)$.
\end{de}

We can associate a canonical {\fs} on $T^{\ssc\cn^{\,\ast}}M$ to each map satisfying {\ca}:
\begin{prop}\label{prop7}
Let $\phi :(M,g)\to (N,h,J)$ be a smooth map from a {\ri} {\mfd} to a {\he} {\mfd}.
\newline If $\phi$ satisfies {\ca} then $(d\phi)^{\ast}(T^{(1,0)\ast}N)$ is an isotropic bundle and therefore defines an {\fs} on $T^{\ssc\cn^{\,\ast}}M$.
\end{prop}

\begin{de}\label{Vde6}
Let $(M,g)$ be a {\ri} {\mfd} and $(N,h,J)$ a {\he} {\mfd}.
\newline Let $\phi :(M,g)\to (N,h,J)$ be a smooth map satisfying {\ca}. We shall call the {\fs} defined by Proposition~\ref{prop7}, the {\em associated {\fs}} and denote it by $F^{\phi}$.
\end{de}

\begin{rk} \label{rk6}
Let $(M,g)$, $(N,h,J)$ and $\phi$ be as in Definition~\ref{Vde6}. Then $\phi :(M,F^{\phi},g)\to (N,h,J)$ is f-{\ho}.
\newline Indeed
\begin{align*}
&(d\phi)^{\ast}(T^{(1,0)\ast}N) \subseteq T^{\ast +}M \\
&\Leftrightarrow \left\{ \begin{array}{c}         
           g^{\ast}\left((d\phi)^{\ast}(T^{(1,0)\ast}N),T^{\ast +}M \right) = \, 0 \\
            g^{\ast}\left((d\phi)^{\ast}(T^{(1,0)\ast}N),T^{\ast 0}M \right) = \, 0 
                               \end{array} \right. \\
&\Leftrightarrow \left\{ \begin{array}{c}
                        (d\phi)(T^{+}M) \subseteq T^{(1,0)}N\\
                        (d\phi)(T^{0}M) \subseteq T^{(0,1)}N 
                         \end{array} \right. 
\end{align*}
but as $(d\phi)(T^{0}M)$ is real, the last equation implies
$$(d\phi)(T^{0}M) = \, 0 .$$
Viewing $J$ as an {\fs}, $(d\phi)^{\ast}(T^{(1,0)\ast}N) \subseteq T^{\ast +}M$ (and its complex conjugate) is equivalent to 
\begin{align*}
&(d\phi)(T^{+}M) \subseteq T^{+}N = T^{(1,0)}N ,\\
&(d\phi)(T^{-}M) \subseteq T^{-}N = T^{(0,1)}N ,\\
&(d\phi)(T^{0}M) \subseteq T^{0}N = \left\{ 0 \right\} ,
\end{align*}
i.e.
$$ d\phi \co F^{\phi} = J \co d\phi ,$$
which means that $\phi$ is f-{\ho}.
\end{rk}

Let $\phi :(M,g)\to (N,h)$ be a smooth map from a {\ri} {\mfd} to a {\he} {\mfd}. Consider $(d\phi)^{\ast}:(\tm ,g^{\ast})\to~(\tn,h^{\ast})$ where $g^{\ast}$ and $h^{\ast}$ are the metrics dual to $g$ and $h$.
\newline Then $(d\phi)^{\ast}$ is a smooth section of the bundle $\left(\phi{^{-1}}\tn\right)^{\ast} \otimes \tm$ which can be identified with the bundle $\tm \otimes \phi{^{-1}}TN$. Let $\nabla$ denote the connection on the bundle $\tm \otimes \phi{^{-1}}TN$ induced from the Levi-Civita connection on $M$ and $N$  then 
$$ \nabla (d\phi)^{\ast}\in C{^{\infty}}\left( \tm \otimes \tm \otimes \phi{^{-1}}TN \right) .$$
Let $X\in TM$ and $Y\in \phi{^{-1}}\tn$. Then
\begin{align*}
\left(\nabla (d\phi)^{\ast}\right) (X,Y) =& \left(\nabla_{X} (d\phi)^{\ast}\right) (Y) \\
=& \nabla^{\stm}_{X} \left( (d\phi)^{\ast}(Y) \right) - (d\phi)^{\ast}\left( \nabla^{\stn}_{d\phi (X)} Y \right),
\end{align*}
which is in $\tm$.
\newline If we denote by $\left(x^{i}\right)_{i=1,\dots,m}$ a system of local coordinates on $M$ in the neighbourhood of the point $x$ and by $\left(z^{\alpha}\right)_{\alpha=1,\dots,n}$ a system of local complex coordinates on $N$ around  the point $\phi(x)$, then a simple computation shows:

\begin{lem}
$$ \sum_{i,j=1}^{m} g^{ij}  \left(\nabla_{\frac{\partial}{\partial x^{i}}} (d\phi)^{\ast}\left(d z^{\alpha} \right) \right) \, \left(\frac{\partial}{\partial x^{j}}\right) = \tau^{\alpha} (\phi), $$
where $\tau(\phi)$ is the tension field.
\end{lem}
\begin{proof}
This follows from the calculation:
\begin{align*}
\nabla_{\frac{\partial}{\partial x^{i}}} (d\phi)^{\ast}\left(d z^{\alpha} \right)
 &= \nabla^{\stm}_{\frac{\partial}{\partial x^{i}}} (d\phi)^{\ast}(d z^{\alpha}) - (d\phi)^{\ast}\left( \nabla^{\stn}_{d\phi (\frac{\partial}{\partial x^{i}})} d z^{\alpha} \right) \\
=& \frac{\partial}{\partial x^{i}}\left( \frac{\partial\phi^{\alpha}}{\partial x^{j}} \right) d x^{j} + \frac{\partial\phi^{\alpha}}{\partial x^{j}} 
\nabla^{\stm}_{\frac{\partial}{\partial x^{i}}} d x^{j} 
- \frac{\partial\phi^{B}}{\partial x^{i}} (d\phi)^{\ast}\left( - {}^N\Gamma^{\alpha}_{BC} dz^{C} \right) \\
=& \left[ \frac{\partial{}^{2}\phi^{\alpha}}{\partial x^{i}\partial x^{j}}  
- {}^{\ssc M}\Gamma^{k}_{ij} \frac{\partial\phi^{\alpha}}{\partial x^{k}}  
+ {}^N\Gamma^{\alpha}_{BC} \frac{\partial\phi^{B}}{\partial x^{i}} 
\frac{\partial\phi^{C}}{\partial x^{j}} \right] dx^{j} .
\end{align*}
\end{proof}

\begin{thm} \label{Vthmpara}
Let $(M,g)$ be a {\ri} {\mfd} and $(N,h,J)$ a {\ka} {\mfd}.
\newline Let $\phi :(M,g)\to (N,h,J)$ be a non-constant smooth map satisfying {\ca}. If the {\fs} $F^{\phi}$ associated to $\phi$ is parallel then $\phi$ is {\ha}.
\end{thm}

\begin{proof}
Using local coordinates $\left(x^{i}\right)_{i=1,\ldots,m}$ on $M$ adapted to the {\fs} $F^{\phi}$, the tension field of $\phi$ can be written:
\begin{align*}
\tau^{\alpha}(\phi) =& \sum_{i,j=1}^{m} g^{ij} \left[ \nabla^{\stm}_{\frac{\partial}{\partial x^{i}}} (d\phi)^{\ast} (d z^{\alpha}) - (d\phi)^{\ast} \left( \nabla^{\stn}_{d\phi (\frac{\partial}{\partial x^{i}})} d z^{\alpha} \right)\right] \left( \frac{\partial}{\partial x^{j}} \right).
\end{align*}
As the tension field $\tau(\phi)$ is real, we shall only consider $\tau^{\alpha}(\phi), \, \alpha \in \{1,\dots,n\}$.

\begin{lem} \label{lem9}
Let $(M,g)$, $(N,h,J)$ and $\phi$ be as in Theorem~\ref{Vthmpara}.
$$\sum_{i,j=1}^{m} g^{ij} (d\phi)^{\ast} \left( \nabla^{\stn}_{d\phi (\frac{\partial}{\partial x^{i}})} d z^{\alpha} \right) \left( \frac{\partial}{\partial x^{j}} \right) = 0 . $$
\end{lem}

\begin{proof}
We only have to consider three cases:
\begin{enumerate}
\item[1)] $\frac{\partial}{\partial x^{i}} \in T^{+}M$ and $\frac{\partial}{\partial x^{j}} \in T^{-}M$,
\item[2)] $\frac{\partial}{\partial x^{i}} \in T^{-}M$ and $\frac{\partial}{\partial x^{j}} \in T^{+}M$,
\item[3)] $\frac{\partial}{\partial x^{i}} , \frac{\partial}{\partial x^{j}} \in T^{0}M$.
\end{enumerate}
1) When $\frac{\partial}{\partial x^{i}} \in T^{+}M$ then $d\phi\left(\frac{\partial}{\partial x^{i}}\right) \in T^{(1,0)}N$ and, since $N$ is {\ka}, 
$$ \nabla^{\stn}_{d\phi (\frac{\partial}{\partial x^{i}})} d z^{\alpha} \in T^{(1,0)\ast}N .$$
Together with Remark~\ref{rk6}, this proves that
$$(d\phi)^{\ast} \left( \nabla^{\stn}_{d\phi (\frac{\partial}{\partial x^{i}})} d z^{\alpha} \right) \in T^{\ast +}M ,$$
but as $\frac{\partial}{\partial x^{j}} \in T^{-}M$,
$$(d\phi)^{\ast} \left( \nabla^{\stn}_{d\phi (\frac{\partial}{\partial x^{i}})} d z^{\alpha} \right) \left( \frac{\partial}{\partial x^{j}} \right) = 0 . $$
2) If $\frac{\partial}{\partial x^{i}} \in T^{-}M$ then $d\phi\left(\frac{\partial}{\partial x^{i}}\right) \in T^{(0,1)}N$ and, since $N$ is {\ka},
$$ \nabla^{\stn}_{d\phi (\frac{\partial}{\partial x^{i}})} d z^{\alpha} = 0 .$$
3) For $\frac{\partial}{\partial x^{i}}\in T^{0}M$, because $\phi$ is f-{\ho}, the vector $d\phi\left(\frac{\partial}{\partial x^{i}}\right)$ vanishes.
\end{proof}

On the other hand, the condition that $F^{\phi}$ be parallel is equivalent to (cf.~\cite[Lemma 2.3]{Raw85}):
$$ \nabla^{\ssc{TM}}_{X} C^{\infty}(T^{+}M) \subseteq C^{\infty}(T^{+}M) \quad\forall X\in TM$$
but, by considering the Christoffel symbols, it is easy to see that:
\begin{align*}
&\nabla^{\ssc{TM}}_{T^{-}M} C^{\infty}(T^{+}M) \subseteq C^{\infty}(T^{+}M) & \Rightarrow & 
\left\{
\begin{array}{l}
\, \nabla^{\ssc{T^{\ast}M}}_{T^{+}M} C^{\infty}(T^{\ast +}M) \subseteq C^{\infty}(T^{\ast +}M + T^{\ast 0}M) \\
\\
 \, \nabla^{\ssc{T^{\ast}M}}_{T^{-}M} C^{\infty}(T^{\ast +}M) \subseteq C^{\infty}(T^{\ast -}M + T^{\ast 0}M)
\end{array} \right. \\
&\nabla^{\ssc{TM}}_{T^{0}M} C^{\infty}(T^{0}M) \subseteq C^{\infty}(T^{0}M) & \Rightarrow & 
\qquad \nabla^{\ssc{T^{\ast}M}}_{T^{0}M} C^{\infty}(T^{\ast +}M) \subseteq C^{\infty}(T^{\ast +}M + T^{\ast -}M).
\end{align*}
Therefore:
\begin{equation}
\nabla^{\stm}_{\frac{\partial}{\partial x^{i}}} (d\phi)^{\ast} (d z^{\alpha}) \in \left\{
\begin{array}{c}
T^{\ast +}M + T^{\ast 0}M  \quad \text{ when } \frac{\partial}{\partial x^{i}} \in T^{+}M, \\
T^{\ast -}M + T^{\ast 0}M \quad \text{ when } \frac{\partial}{\partial x^{i}} \in T^{-}M, \\
T^{\ast +}M + T^{\ast -}M \quad \text{ when } \frac{\partial}{\partial x^{i}} \in T^{0}M.
\end{array}
\right.
\end{equation}
Since $g^{ij} =0$ when $\frac{\partial}{\partial x^{i}}$ and $\frac{\partial}{\partial x^{j}}$ are both of the same type $+$ or $-$, or when one is of type $0$ and the other is of type $+$ or $-$, we have that
\begin{align*}
\tau^{\alpha}(\phi) =& \sum_{i,j=1}^{m} g^{ij} \left[ \nabla^{\stm}_{\frac{\partial}{\partial x^{i}}} (d\phi)^{\ast} (d z^{\alpha}) \right] \left( \frac{\partial}{\partial x^{j}} \right) = 0.
\end{align*}
\end{proof}

\begin{rk}
The condition on $F^{\phi}$ in Theorem~\ref{Vthmpara} is a strong one. It is the analogue of being {\ka} for an almost {\he} {\mfd}.
\end{rk}

We can weaken the condition on the associated {\fs} if we introduce a condition on the fundamental 2-form $\omega$.

\begin{de}~\cite{KonYano84}
Let $(M^{m},g)$ be a {\ri} {\mfd} with an {\fs} $F$. We define the {\em fundamental 2-form} of $(M^{m},g,F)$ to be,
$$\omega(X,Y) = g(X,FY) \qquad \forall X,Y \in T^{\ssc\cn}M.$$
In a frame $(\theta_{i})_{i=1,\ldots,m}$ adapted to the {\fs} $F$:
$$\omega = \frac{1}{2}\sum_{\substack{i<j = 1 \\ i,j \in (+)}}^{m} 
g( \theta_{i},\theta_{\overline{\jmath}}) \, \Theta^{i}\wedge \Theta^{\overline{\jmath}},$$
where $i \in (+)$ means that $\theta_{i} \in T^{+}M$, $\left\{\Theta^{i}\right\}$ is the dual frame of $\left\{\theta_{i}\right\}$ and, for $i \in (+)$, $\Theta^{\overline{\imath}}= \overline{\Theta^{i}} \in T^{\ast -}M$.
\end{de}

\begin{thm} \label{prop}
Let $(M,g)$ be a {\ri} {\mfd} and $(N,h,J)$ a {\ka} {\mfd}.
\newline Let $\phi :(M,g)\to (N,h,J)$ be a non-constant smooth map satisfying {\ca}.
Suppose that the associated {\fs} $F^{\phi}$ is integrable, and
\begin{equation}\label{met}
\nabla_{T^{0}M} \left( T^{\ast +}M \right) \subseteq T^{\ast +}M + T^{\ast -}M,
\end{equation}
and the fundamental 2-form $\omega$ is such that:
\begin{equation}
d\omega ^{(1,2)}  = 0, \label{sym}
\end{equation}
where $d\omega ^{(1,2)} = 0$ means that $d\omega (u,v,w) = 0$ whenever $u,v$ are of the same type and $w$ of a different type.
\newline Then $\phi$ is a {\ha} map.
\end{thm}

\begin{proof}
Let $(x^{j})$, $(z^{\alpha})$ be adapted local coordinates on $M$ and $N$.
Recall that the tension field of $\phi$ can be written as:
$$\tau^{\alpha}(\phi) =\sum_{i,j=1}^{m} g^{ij} \left[ \nabla^{\stm}_{\frac{\partial}{\partial x^{i}}} (d\phi)^{\ast} (d z^{\alpha}) - (d\phi)^{\ast} \left( \nabla^{\stn}_{d\phi (\frac{\partial}{\partial x^{i}})} d z^{\alpha} \right)\right] \left( \frac{\partial}{\partial x^{j}} \right)$$
and that (Lemma~\ref{lem9}):
$$(d\phi)^{\ast} \left( \nabla^{\stn}_{d\phi (\frac{\partial}{\partial x^{i}})} d z^{\alpha} \right) = 0 . $$
It is clear that Condition~\eqref{met} implies
\begin{equation*}
\nabla^{\stm}_{\frac{\partial}{\partial x^{i}}} (d\phi)^{\ast} (d z^{\alpha}) \in \, T^{\ast +}M + T^{\ast -}M  \quad \forall i \, \text{such that} \quad \frac{\partial}{\partial x^{i}} \in T^{0}M .
\end{equation*}
Note that Condition~\eqref{met} is equivalent to
$$\frac{\partial g_{ij}}{\partial x^{l}} = 0 \quad \forall i, j \in (0), l \in (-).$$
It only remains to show that Condition~\eqref{sym} implies: 
\begin{subequations}\label{VA}
\begin{align}
&\nabla^{\stm}_{\frac{\partial}{\partial x^{i}}} (d\phi)^{\ast} (d z^{\alpha}) \in T^{\ast +}M + T^{\ast 0}M \quad \forall i \, \text{such that} \quad \frac{\partial}{\partial x^{i}} \in T^{+}M ,\\
&\nabla^{\stm}_{\frac{\partial}{\partial x^{i}}} (d\phi)^{\ast} (d z^{\alpha}) \in T^{\ast -}M + T^{\ast 0}M \quad \forall i \, \text{such that} \quad\frac{\partial}{\partial x^{i}} \in T^{-}M .
\end{align}
\end{subequations}

\begin{lem}  \label{lem10}
The fundamental 2-form $\omega$ satisfies:
\begin{equation*}
d\omega ^{(1,2)} = 0
\end{equation*}
if and only if
\begin{equation}\label{equa1}
\nabla^{\stm}_{v} C^{\infty}(T^{\ast +}M + T^{\ast -}M) \in C^{\infty}(T^{\ast +}M + T^{\ast 0}M) \qquad \forall v\in C^{\infty}(T^{+}M).
\end{equation}
\end{lem}

\begin{proof}
Take adapted coordinates $(x^{j})$.
Suppose, throughout the proof, that $\frac{\partial}{\partial x^{j}} \in T^{+}M$.
\newline It is easy to see that~\eqref{equa1} is equivalent to ${}^{\ssc M}\Gamma^{j}_{ik} = 0$ whenever $dx^{i}$ has type $+$ or $-$ and $dx^{k}$ has type $-$ or $+$, different from $dx^{i}$.
\newline We can express this condition in terms of $\omega$.
\newline Let us first compute $d\omega$:
\begin{align*}
d\omega =& \sum_{k=1}^{m} \sum_{\substack{i<j = 1 \\ i,j \in (+)}}^{m} \frac{\partial g_{i\overline{\jmath}}}{\partial x^{k}}  \, dx^{k}\wedge dx^{i}\wedge dx^{\overline{\jmath}} \\
=& \sum_{\substack{k<i<j = 1 \\ k,i,j \in (+) }}^{m} \left( \frac{\partial g_{i\overline{\jmath}}}{\partial x^{k}} - \frac{\partial g_{k\overline{\jmath}}}{\partial x^{i}} + \frac{\partial g_{ik}}{\partial x^{\overline{\jmath}}} \right) \, dx^{k}\wedge dx^{i}\wedge dx^{\overline{\jmath}} \\
&+ \sum_{\substack{k,i<j = 1 \\ i,j \in (+), k\in (0)}}^{m}  \frac{\partial g_{i\overline{\jmath}}}{\partial x^{k}}  \, dx^{k}\wedge dx^{i}\wedge dx^{\overline{\jmath}} .
\end{align*}

\begin{notation}
$dx^{k_{\alpha}} \in T^{\alpha}M$ for $\alpha= +, -, 0$.
\end{notation}
It is easy to see that $d\omega^{(1,2)}=0$ is equivalent to
\begin{align*}
\frac{\partial g_{i_{+}l_{-}}}{\partial x^{k_{-}}} = \frac{\partial g_{i_{+}k_{-}}}{\partial x^{l_{-}}} .
\end{align*}
Since
\begin{align*}
{}^{\ssc M}\Gamma^{j}_{ik} =& \frac{1}{2} 
g^{jl_{-}} 
\left(
  \frac{\partial g_{l_{-}k}}{\partial x^{i}} 
- \frac{\partial g_{ik}}{\partial x^{l_{-}}} 
+ \frac{\partial g_{il_{-}}}{\partial x^{k}}
\right) ,
\end{align*}
we first observe that, for $dx^{i}\in T^{-}M$ and $dx^{k}\in T^{0}M$,
$${}^{\ssc M}\Gamma^{j}_{ik} =0,$$
and that
$d\omega^{(1,2)}=0$ implies
$${}^{\ssc M}\Gamma^{j}_{ik} =0$$
in both cases
\begin{enumerate}
\item $dx^{i}\in T^{+}M$ and $dx^{k}\in T^{-}M$,
\item $dx^{i}\in T^{-}M$ and $dx^{k}\in T^{+}M$.
\end{enumerate}
Hence Equation~\eqref{sym} implies that:
\begin{align*}
&\nabla_{T^{+}M} \left( T^{\ast +}M \right) \subseteq T^{\ast +}M + T^{\ast 0}M ,\\
&\nabla_{T^{-}M} \left( T^{\ast +}M \right) \subseteq T^{\ast -}M + T^{\ast 0}M.
\end{align*}
\end{proof}

To end the proof of Theorem~\ref{prop}, we deduce from Lemma~\ref{lem10} that \eqref{VA} holds. But $g^{ij} = 0$ for $i$ and $j$ both of type $+$ or $-$ or if $i$ is of type $0$ and $j$ is not. Thus
\begin{equation*}
\tau^{\alpha}(\phi) = \sum_{i,j=1}^{m} 
g^{ij}
 \left[
 \nabla^{\stm}_{\frac{\partial}{\partial x^{i}}} (d\phi)^{\ast}(d z^{\alpha}) \right] 
\left(
\frac{\partial}{\partial x^{j}} 
\right) = 0
\end{equation*}
and, as $\tau(\phi)$ is real,

$$\tau(\phi)\equiv 0.$$

\end{proof}

\begin{cor}
Let $(M,g)$ carry an integrable {\fs} satisfying Equation~\eqref{met} and whose fundamental 2-form satisfies~\eqref{sym}, and let $(N,J,h)$ be an {\ka} {\mfd}. Then any f-{\ho} map $\phi: M \to N$ is {\ha}.
\end{cor}


\begin{thebibliography}{BBdBR89}

\bibitem[BBdBR89]{BurBurdeB89}
D.~Burns, F.E. Burstall, P.~de~Bartolomeis, and J.~Rawnsley.
\newblock {S}tability of harmonic maps of {K}\"ahler manifolds.
\newblock {\em J. Differential Geom.}, 30:579--594, 1989.

\bibitem[Che93]{Che93A}
J.Y. Chen.
\newblock {S}table harmonic maps into {${S}^2$}.
\newblock In T.~Kotake, S.~Nishikawa, and R.~Schoen, editors, {\em {G}eometry
  and {G}lobal {A}nalysis}, pages 431--436. T\^ohoku Univ., Sendai, 1993.

\bibitem[IY64]{IshiharaYano64}
S.~Ishihara and K.~Yano.
\newblock On integrability conditions of a structure $f$ satisfying $f^{3} + f
  = 0$.
\newblock {\em Quart. J. Math. Oxford}, 15(2):217--222, 1964.

\bibitem[KY84]{KonYano84}
M.~Kon and K.~Yano.
\newblock {\em Structures on Manifolds}, volume~3 of {\em Series in Pure
  Mathematics}.
\newblock World Scientific, 1984.

\bibitem[OU90]{OhnUda90}
Y.~Ohnita and S.~Udagawa.
\newblock {S}tability, complex-analyticity, and constancy of pluriharmonic maps
  from compact {K}\"ahler manifolds.
\newblock {\em Math.\ Z.}, 205:624--644, 1990.

\bibitem[OV90]{OhnVal90}
Y.~Ohnita and G.~Valli.
\newblock {P}luriharmonic maps into compact {L}ie groups and factorization into
  unitons.
\newblock {\em Proc.\ London Math.\ Soc.}, 61:546--570, 1990.

\bibitem[Raw85]{Raw85}
J.~Rawnsley.
\newblock {$f$}-structures, {$f$}-twistor spaces and harmonic maps.
\newblock In E.~Vesentini, editor, {\em {G}eometry {S}eminar ``{L}uigi
  {B}ianchi'' {I}{I}---1984}, volume 1164 of {\em Lecture Notes in Math.},
  pages 85--159. Springer, Berlin, Heidelberg, New York, 1985.

\bibitem[Uda88]{Uda88A}
S.~Udagawa.
\newblock {P}luriharmonic maps and minimal immersions of {K}{\"a}hler
  manifolds.
\newblock {\em J. London Math.\ Soc.}, 37:375--384, 1988.

\bibitem[Yan63]{Yano63}
K.~Yano.
\newblock On a structure defined by a tensor field of type {$(1,1)$} satisfying
  {$f^{3} + f = 0$}.
\newblock {\em Tensor N.S.}, 14:99--109, 1963.

\end{thebibliography}
\end{document}